\documentclass[11pt]{amsart}

\usepackage[T1]{fontenc}
\usepackage{lmodern}
\usepackage[margin=1in]{geometry}
\usepackage{amsmath,amssymb,amsthm,mathtools}
\usepackage{letltxmacro}
\usepackage[hidelinks]{hyperref}
\usepackage{microtype}

\usepackage{newtx}
\usepackage{xcolor}

\newcommand{\R}{\mathbb{R}}
\newcommand{\C}{\mathbb{C}}
\newcommand{\Sn}{\mathbb{S}^{n}}
\newcommand{\Tn}{\mathbb{T}^{n}}
\newcommand{\ip}[2]{\left\langle #1,#2\right\rangle}
\DeclareMathOperator{\divg}{div}

\newtheoremstyle{manuscript}
  {6pt}{6pt}{\itshape}{}{\bfseries}{.}{0.5em}{}
\theoremstyle{manuscript}
\newtheorem{theorem}{Theorem}

\newtheorem{conjecture}{Conjecture}
\newtheorem{corollary}{Corollary}
\newtheorem{proposition}{Proposition}
\newtheorem{lemma}{Lemma}
\theoremstyle{remark}
\newtheorem*{remark}{{\bf Remark}}

\newenvironment{resultparts}{%
  \leavevmode\par\nobreak
  \begin{list}{\textup{(\roman{enumi})}}{%
    \usecounter{enumi}%
    \setlength{\leftmargin}{2em}%
    \setlength{\labelwidth}{1.5em}%
    \setlength{\labelsep}{0.5em}%
    \setlength{\topsep}{0pt}%
    \setlength{\partopsep}{0pt}%
    \setlength{\itemsep}{2pt}%
    \setlength{\parsep}{0pt}%
  }%
}{\end{list}}

\LetLtxMacro{\originalproof}{\proof}
\let\endoriginalproof\endproof
\renewenvironment{proof}[1][\proofname]{%
  \originalproof[\normalfont\bfseries #1]%
}{\endoriginalproof}
\newcommand{\proofpart}[1]{\noindent\textup{(#1)}\enspace}

\title[The Whitney sphere and tori]{The Whitney sphere and tori as self-shrinkers}
\author{Jingyi Chen}
\address{Department of Mathematics, The University of British Columbia, Vancouver, BC, Canada}
\email{jychen@math.ubc.ca}
\author{Micah Warren}
\address{Department of Mathematics, University of Oregon, Eugene, OR, USA}
\email{micahw@uoregon.edu}

\thanks{J. Chen is supported by NSERC Discovery Grant GR010074.}
\date{}

\begin{document}

\begin{abstract}
We prove that the Whitney sphere in $\C^n$ is a self-shrinker for
the Hamiltonian gradient flow of volume. We construct an explicit
family of Lagrangian self-shrinking tori in $\C^2$ containing infinitely
many pairwise noncongruent examples. We also characterize closed
connected Lagrangian immersions whose Legendrian lifts lie in a
Kor\'anyi sphere by showing that, up to scaling and a unitary
transformation, they factor through the Whitney immersion via a
diffeomorphism for $n\geq2$ and a finite covering for $n=1$.
\end{abstract}

\maketitle

\section{Introduction}

In this note, we prove that the Whitney sphere is a self-shrinker for the
Hamiltonian gradient flow of volume. We also construct explicit
self-shrinking tori and prove a rigidity theorem characterizing the
Whitney immersion. The flow, studied recently by the authors in
\cite{ChenWarren}, is given by
\begin{equation}\label{eq:flow}
\partial_tF=J\nabla_g\divg_g(JH).
\end{equation}
We work in $\C^n$ with the standard complex structure $J$.
All differential operators are taken in the induced metric $g$, with
volume form $dV_g$, and $\Delta_g=\divg_g\nabla_g$.
We denote the tangential and normal components of a vector along the
immersion by $V^\top$ and $V^\perp$.
For Lagrangian immersions $F:M\times[0,T)\to\C^n$ from a closed $n$-dimensional manifold $M$, \eqref{eq:flow} satisfies
\[
\frac{d}{dt}\operatorname{Vol}(F)
=-\int_M\bigl|\divg_g(JH)\bigr|^2\,dV_g.
\]
This is the negative gradient flow of volume for the $L^2$ metric on
Hamiltonian potentials. Its stationary solutions are the Hamiltonian
stationary immersions of Oh \cite{Oh1993}, characterized by a locally
harmonic Lagrangian angle \cite{SchoenWolfson1999}.
For a locally defined Lagrangian angle $\theta$, the mean curvature is
$H=J\nabla_g\theta$, so the velocity in \eqref{eq:flow} is
$-J\nabla_g\Delta_g\theta$. Equation~\eqref{eq:flow} is a fourth-order
geometric evolution system for $F$, but is not strictly parabolic.
Locally, writing the Lagrangian as the graph of $Du$, it reduces, up to
reparametrization and normalization of $u$, to the fourth-order strictly
parabolic scalar equation $u_t=-\Delta_g\theta$; see
\cite[Section~3.1]{ChenWarren}.
%For immersions, the Hamiltonian potentials are functions on the domain; see \cite[Section 2.2]{ChenWarren}. 
%When $N=\C^n$, for a locally defined Lagrangian angle $\theta$, $H=J\nabla_g\theta$ by \cite{Harvey-Lawson}, so the velocity is $-J\nabla_g\Delta_g\theta$.

A Lagrangian immersion $F$ is a \emph{self-shrinker} of \eqref{eq:flow},
centred at the origin, if
\[
J\nabla_g\divg_g(JH)=-cF^\perp
\]
for some $c>0$. Scaling and reparametrization then give a shrinking
solution of \eqref{eq:flow}.

Smoczyk proved that closed Lagrangian self-shrinkers for mean
curvature flow have nontrivial Maslov class
\cite[Theorem 2.3.5 and Corollary 2.3.6]{Smoczyk2000}; see also
\cite[p.~265]{Smoczyk2012}; in particular, 
the Whitney sphere is not self-shrinking under mean curvature flow, since $H^1(\Sn;\R)=0$ for $n>1$.
For $n=1$, the Bernoulli lemniscate is not a self-shrinker for
curve shortening flow, by the classification of Abresch and Langer
\cite[Theorem~A]{AbreschLanger1986}.

For $n\geq1$, write
$\Sn=\{(u,s)\in\R^n\times\R:|u|^2+s^2=1\}$ and define the Whitney immersion
\begin{equation}\label{eq:whitney}
W_n(u,s)=\frac{1+is}{1+s^2}\,u.
\end{equation}
We write $W=W_n$ when no confusion arises.

\begin{theorem}\label{thm:main}
For every $n\geq1$, the Whitney immersion \eqref{eq:whitney} is a
self-shrinker for the flow \eqref{eq:flow}, with
\begin{equation}\label{eq:soliton}
J\nabla_g\divg_g(JH)=-2n(n+2)W^\perp.
\end{equation}
\end{theorem}

The flow splits over Cartesian products, giving the following consequence.

\begin{corollary}[Products of Whitney spheres]\label{cor:product-spheres}
Let $n_i\geq1$ and $a_i>0$ for $1\leq i\leq k$, and set
$M=\prod_{i=1}^k\mathbb{S}^{n_i}$ and $N=\sum_{i=1}^k n_i$.
For $(u_i,s_i)\in\mathbb{S}^{n_i}\subset\R^{n_i}\times\R$,
the product immersion $F_0\colon M\to\C^N$ defined by
\[
F_0\bigl((u_1,s_1),\ldots,(u_k,s_k)\bigr)
=\bigl(a_1W_{n_1}(u_1,s_1),\ldots,a_kW_{n_k}(u_k,s_k)\bigr)
\]
is an exact Lagrangian immersion, and it is not embedded. It is a
self-shrinker for \eqref{eq:flow} if and only if
\begin{equation}\label{eq:product-balance}
\frac{2n_i(n_i+2)}{a_i^4}=c
\quad\text{for every }i
\end{equation}
for a common constant $c>0$, that is, if and only if all factors have the
same collapse time $a_i^4/(8n_i(n_i+2))$.
\end{corollary}

Products of round circles, possibly with different radii, are Hamiltonian
stationary and hence are not self-shrinkers. The following theorem gives
both products of planar curves and tori that are not congruent to such
products. The construction in part~(ii) uses that of Castro and Lerma
\cite{CastroLerma2014}.

\begin{theorem}\label{thm:circle-tori}
For the flow \eqref{eq:flow}, 
\begin{resultparts}
\item
For every $n\geq1$, the Cartesian product of $n$ copies of the
Bernoulli lemniscate $W_1$ is an exact, immersed self-shrinker $\Tn\to\C^n$.
\item
For every integer $N\geq1$ and coprime positive integers $m_1,m_2$,
the explicit torus $F_{N;m_1,m_2}$ defined in
\eqref{eq:nonproduct-curve}--\eqref{eq:weighted-torus} is an exact
self-shrinker in $\C^2$, with constant
$c_N=4/[N(N+1)]$. The torus has nonflat induced metric and is not
congruent to a Cartesian product of planar curves.
Tori with different $N$ are not related by unitary transformations,
dilations and reparametrizations. For each $N$, the subfamily
$F_{N;1,m}$, $m\geq1$, is pairwise noncongruent.
\item
Let $z\colon I\to\C\setminus\{0\}$ be a smooth regular curve,
where $I$ is an interval or a circle, and let $\psi$ be a smooth
closed regular Legendrian curve in $\mathbb{S}^3\subset\C^2$.
If $F(\sigma,\tau)=z(\sigma)\psi(\tau)$ is a self-shrinker for \eqref{eq:flow}, then 
its image is invariant under a nontrivial circle subgroup of $U(2)$.
\end{resultparts}
\end{theorem}

The shrinker criterion in Lemma~\ref{lem:criterion} implies exactness.
Gromov's theorem \cite{Gromov1985} then rules out closed embedded
self-shrinkers.
\begin{theorem}\label{thm:exactness}
Every self-shrinker for \eqref{eq:flow} is exact.
In particular, there are no closed embedded self-shrinkers
in $\C^n$.
\end{theorem}

By Theorem~\ref{thm:exactness}, every self-shrinker admits a Legendrian
lift to the Heisenberg group $\C^n\times\R$. Such a lift is a pair
$(F,\varphi)$ satisfying $d\varphi=\ip{JF}{dF}$; its Kor\'anyi norm is
$\mathfrak r=\bigl(|F|^4+4\varphi^2\bigr)^{1/4}$. 
The next theorem characterizes the Whitney immersion by requiring its
Legendrian lift to lie in a Kor\'anyi sphere.

\begin{theorem}[Rigidity in a Kor\'anyi sphere]\label{thm:koranyi-rigidity}
Let $F\colon M^n\to\C^n$, $n\geq1$, be a smooth Lagrangian immersion
of a closed connected manifold. Suppose $F$ admits a Legendrian
lift $(F,\varphi)$ contained in the Kor\'anyi sphere of radius $R>0$
centered at the origin,
\[
S_K(R)=\{(z,t)\in\C^n\times\R:|z|^4+4t^2=R^4\}.
\]
For $n\geq2$, there exist $U\in U(n)$ and a diffeomorphism
$\Pi\colon M\to\mathbb{S}^n$ such that
\[
F=R\,U\circ W_n\circ\Pi.
\]
For $n=1$, the same factorization holds with $U\in U(1)$ and a
finite smooth covering $\Pi\colon M\to\mathbb{S}^1$. Thus $F$ is a
finite cover of a scaled and rotated Bernoulli lemniscate.

In both cases, $F$ is a self-shrinker for \eqref{eq:flow}
with constant $c=2n(n+2)/R^4$.
\end{theorem}

The Whitney sphere and the tori above provide examples in genus zero
and one. We conjecture that such examples exist in every genus.

\begin{conjecture}
For every integer $g\geq0$, there exists a closed connected orientable
surface $M_g$ of genus $g$ and a Lagrangian immersion
$F\colon M_g\to\C^2$ that is a self-shrinker for the flow \eqref{eq:flow}.
\end{conjecture}
Section~\ref{subsec:heisenberg} derives integral identities for closed
shrinkers from their normalized lifts, proves
Theorem~\ref{thm:koranyi-rigidity}, and describes the lifts of the
Whitney sphere and the tori $F_{N;m_1,m_2}$.
The coordinates used there also appear in the monotonicity formulas
for Hamiltonian stationary Legendrian surfaces of Schoen and Wolfson
\cite{SchoenWolfson1999,SchoenWolfson2001} and Pigati and Rivi\`ere
\cite{PigatiRiviere}.

\smallskip

\noindent\textbf{AI Usage Statement:} The authors formulated the research
questions and directed the development of this work, building on their
earlier work \cite{ChenWarren}. AI tools assisted with mathematical
exploration, computations, proof development, and manuscript revision.
The authors have verified the proofs and take full responsibility for the
content and conclusions of this paper.

\smallskip

\noindent\textbf{Formal verification:} All results of this paper are
verified in the Lean~4 proof assistant using the Mathlib library, with
no axioms beyond those of Lean's core logic, except for the use of
Gromov's theorem \cite[Corollary 2.3.B$_2$]{Gromov1985} as an explicit
hypothesis in the second sentence of Theorem~\ref{thm:exactness} and the
exclusion of an embedded endpoint in
Proposition~\ref{prop:deformation-obstructions}(ii).
Gromov's theorem itself is not formalized. The verification includes the
global results Theorem~\ref{thm:koranyi-rigidity} and
Proposition~\ref{prop:heisenberg}, and, through a formalization of
\cite[Theorem 3.4]{ChenWarren}, the homothetic evolution and its
uniqueness for the Whitney sphere, for products of two Whitney spheres
with equal collapse times, and for the tori $F_{N;m_1,m_2}$. The code is
available at \texttt{https://github.com/micah541/HSLflow-Whitney}.

\section{Proofs and geometric consequences}\label{sec:proof}

Throughout, the domain manifold $M$ is assumed connected.
On $\C^n=\R_q^n\times\R_p^n$, we use the Euclidean metric, the complex structure 
$J(q,p)=(-p,q)$, and the symplectic structure 
\[
\omega_0=\sum_{j=1}^n dq_j\wedge dp_j,\quad
\lambda_0=\frac12\sum_{j=1}^n(q_j\,dp_j-p_j\,dq_j),
\quad d\lambda_0=\omega_0.
\]
A Lagrangian immersion $F$ is exact if $F^*\lambda_0$ is an exact
one-form on its domain. 

%All the shrinker equations below are checked with the following criterion.

\begin{lemma}[Shrinker criterion]\label{lem:criterion}
Let $F\colon M\to\C^n$ be a Lagrangian immersion and $c\in\R$. Then
\[
J\nabla_g\divg_g(JH)=-cF^\perp
\quad\Longleftrightarrow\quad
d\,\divg_g(JH)=2c\,F^*\lambda_0 .
\]
If $F^*\lambda_0=df$ and $\theta$ is a local Lagrangian angle,
these conditions are equivalent to $\Delta_g\theta+2cf$ being locally constant.
\end{lemma}

\begin{proof}[Proof of Lemma~\ref{lem:criterion}]
For $z\in\C^n$ and $v\in T_z\C^n$, we have $\lambda_0(v)=\frac12\ip{Jz}{v}$.
Since $F$ is Lagrangian, $J$ maps the tangent space $dF(TM)$ isometrically
onto the normal space and vice versa, so $J(V^\perp)=(JV)^\top$ for every
$V\in\C^n$. Set $h=\divg_g(JH)$. Applying $J$, the equation
$JdF(\nabla_g h)=-cF^\perp$ is equivalent to $dF(\nabla_g h)=c(JF)^\top$.
Taking inner products with $dF(X)$, this is equivalent to
\[
dh(X)=c\ip{JF}{dF(X)}=2c\,(F^*\lambda_0)(X)
\quad\text{for every tangent vector }X\text{ on }M.
\]
For the last statement, the expression $H=J\nabla_g\theta$ of \cite{HL1982} gives $JH=-\nabla_g\theta$, 
so $h=-\Delta_g\theta$. 
\end{proof}

Even when the Maslov class is nonzero and $\theta$ exists only locally,
$\nabla_g\theta=-JH$ and $\Delta_g\theta=-\divg_g(JH)$ are globally defined.
All identities involving derivatives of $\theta$ are understood in this sense.

\begin{remark}[Legendrian lift]
If $F^*\lambda_0=df$, put $\varphi=2f$. Then $(F,\varphi)$ is a Legendrian
immersion into the Heisenberg group $\C^n\times\R$ with the contact form
$d\varphi-2\lambda_0$. For $n=2$ this is the normalization of
\cite{PigatiRiviere}. For fixed $F$, the height $\varphi$ is unique up to
an additive constant, corresponding to a \emph{vertical translation}
in the $\R$ direction. For a self-shrinker with constant $c>0$,
Lemma~\ref{lem:criterion} uniquely fixes this constant by requiring
\begin{equation}\label{eq:legendrian-shrinker}
\Delta_g\theta=-c\,\varphi ,\qquad\text{equivalently}\qquad
\divg_g(JH)=c\,\varphi .
\end{equation}
\end{remark}

\subsection{The Whitney sphere}\label{subsec:whitney}

\begin{proof}[Proof of Theorem~\ref{thm:main}]
Let $g_0$ be the unit round metric on $\Sn$. We first verify that
the Whitney immersion is Lagrangian and satisfies
\begin{align}
g&=\frac{1}{1+s^2}g_0,\label{eq:metric}\\
\theta&=(n+2)\arctan s+\frac\pi2,\label{eq:angle}\\
W^*\lambda_0&=df_W,\quad
f_W(u,s)=\frac{s}{2(1+s^2)},\label{eq:primitive}\\
\Delta_g\theta&=-4n(n+2)f_W,\label{eq:lap-angle}
\end{align}
where $\theta$ is a globally defined real-valued Lagrangian angle.

By smoothness, it suffices to verify these identities away from the
north pole $(0,1)$, using inverse stereographic projection
\[
y\longmapsto(u,s)=\left(\frac{2y}{1+|y|^2},\,\frac{|y|^2-1}{|y|^2+1}\right),
\qquad y\in\R^n,
\]
a diffeomorphism onto $\Sn\setminus\{(0,1)\}$. Put $t=|y|^2$, so that
$s=(t-1)/(t+1)$. In these coordinates, \eqref{eq:whitney} becomes
\[
W(y)=\psi(t)\,y,\qquad \psi(t)=\frac{1-i}{1-it},\qquad
\psi'=\frac{i\psi}{1-it}.
\]
Let $h(z,w)=\sum_j\overline{z_j}w_j=\ip zw+i\,\omega_0(z,w)$ be the
Hermitian product on $\C^n$. Since $dW(a)=\psi a+2\psi'\ip ya\,y$ for
$a\in\R^n$,
\[
h\bigl(dW(a),dW(b)\bigr)=|\psi|^2\ip ab
+4\bigl(\operatorname{Re}(\overline\psi\psi')+t|\psi'|^2\bigr)\ip ya\ip yb .
\]
This is real, so $W^*\omega_0=0$. Moreover
$\operatorname{Re}(\overline\psi\psi')=-t|\psi|^2/(1+t^2)$ and
$|\psi'|^2=|\psi|^2/(1+t^2)$, so the bracket vanishes and
\[
g=|\psi|^2|dy|^2=\frac{2\,|dy|^2}{1+t^2}.
\]
Since $g_0=4|dy|^2/(1+t)^2$ and $1+s^2=2(1+t^2)/(1+t)^2$, this is
\eqref{eq:metric}.

Since $\lambda_0=\frac12\operatorname{Im}h(z,dz)$ and $\ip y{dy}=\frac12dt$,
\[
W^*\lambda_0
=\frac12\operatorname{Im}\bigl((|\psi|^2+2t\,\overline\psi\psi')\ip y{dy}\bigr)
=\frac{t\,dt}{(1+t^2)^2}
=d\left(\frac{t^2-1}{4(1+t^2)}\right),
\]
and $(t^2-1)/(4(1+t^2))=s/(2(1+s^2))$. This proves \eqref{eq:primitive}.

The columns of $dW$ form the complex matrix $\psi I+2\psi'\,yy^{T}$, so
\[
\det_\C dW=\psi^{n-1}(\psi+2t\psi')=\frac{(1-i)^n(1+it)}{(1-it)^{n+1}}.
\]
For a Lagrangian immersion $|\det_\C dW|=\sqrt{\det g}$, so the definition
$W^*\Omega=e^{i\theta}dV_g$, with $\Omega=\bigwedge_{j=1}^n(dq_j+i\,dp_j)$,
gives
\[
\theta=\arg\det_\C dW=(n+2)\arctan t-\frac{n\pi}{4}.
\]
Since $\arctan t=\arctan s+\pi/4$, this is \eqref{eq:angle}, whose right
side is smooth on all of $\Sn$.

Finally, for a conformal metric $g=\rho\,|dy|^2$ on a domain in $\R^n$,
$\Delta_g=\rho^{-n/2}\sum_i\partial_i(\rho^{n/2-1}\partial_i)$. With
$\rho=2/(1+t^2)$ and $\partial_i\theta=2(n+2)y_i/(1+t^2)$,
\[
\Delta_g\theta
=(n+2)(1+t^2)^{n/2}\sum_{i=1}^n\partial_i\Bigl(\frac{y_i}{(1+t^2)^{n/2}}\Bigr)
=\frac{n(n+2)(1-t^2)}{1+t^2}
=-4n(n+2)f_W.
\]

By \eqref{eq:lap-angle}, $\Delta_g\theta+2cf_W=0$ with $c=2n(n+2)$, so
Lemma~\ref{lem:criterion} gives \eqref{eq:soliton}.
\end{proof}

The north and south poles map by $W$ to the unique transverse double point
at the origin. Their primitive values are
$f_W(0,\pm1)=\pm1/4$, so $f_W$ does not descend to the image.

For a self-shrinker $F_*\colon M\to\C^n$ with $M$ closed and constant
$c>0$, the velocity in \eqref{eq:flow} scales by $r^{-3}$ under dilation
by $r$. Thus $r(t)F_*$ satisfies the normal velocity equation when
\[
r'=-cr^{-3},\qquad r(t)^4=r_0^4-4ct,\qquad
0\leq t<T:=\frac{r_0^4}{4c},
\]
where $r(0)=r_0>0$. To cancel the tangential velocity, let
$dF_*(X)=F_*^\top$ and solve
\[
\partial_t\phi_t=-\frac{r'}{r}X\circ\phi_t,
\qquad \phi_0=\operatorname{id}_M.
\]
Compactness of $M$ gives $\phi_t$ for all $t<T$. Then
$\widetilde F_t=r(t)F_*\circ\phi_t$ satisfies
$\partial_t\widetilde F_t=r'F_*^\perp\circ\phi_t$, hence
\eqref{eq:flow} with initial immersion $r_0F_*$.
See \cite[Claim 3.3 and Theorem 3.4]{ChenWarren} for reparametrization
and uniqueness. Since $\widetilde F_t\to0$ uniformly as $t\uparrow T$,
the solution has no smooth extension through $T$ as an immersion.
Mean curvature scales by $r^{-1}$, so for $H_*\not\equiv0$,
\[
\max_M|H_t|=(4c)^{-1/4}\max_M|H_*|\,(T-t)^{-1/4}.
\]
After the rescaling in \cite[\S~4.1]{Chen2025}, these homothetic
solutions become independent of time up to reparametrization.
This applies to the Whitney sphere, the tori below, and the products
satisfying \eqref{eq:product-balance}.

For $F_*=W$, the scale and collapse time are
\begin{equation}\label{eq:radius}
r(t)=\bigl(r_0^4-8n(n+2)t\bigr)^{1/4},\qquad
0\leq t<T=\frac{r_0^4}{8n(n+2)}.
\end{equation}

The one-dimensional Whitney immersion is the Bernoulli lemniscate
\begin{equation}\label{eq:lemniscate}
\gamma(\alpha):=W(\cos\alpha,\sin\alpha)
=\frac{\cos\alpha(1+i\sin\alpha)}{1+\sin^2\alpha},
\quad \alpha\in\R/2\pi\mathbb{Z}.
\end{equation}

\begin{remark}[The curve diffusion case]
For $\gamma=q+ip$, the Cartesian equation is
\[
(q^2+p^2)^2=q^2-p^2.
\]
The transverse double point at the origin has preimages
$\alpha=\pi/2,3\pi/2$ and tangent lines $p=\pm q$.

Let $\ell$ denote arclength and $\kappa$ the signed curvature
with respect to the unit normal $\nu=J\gamma_\ell$.
Then $H=\kappa\nu$ and
\[
JH=-\kappa\gamma_\ell,\quad
\divg(JH)=-\kappa_\ell,\quad
J\nabla\divg(JH)=-\kappa_{\ell\ell}\nu.
\]
Hence \eqref{eq:flow} reduces to curve diffusion. A round circle has
constant curvature and is stationary under this flow.

For $n=1$, \eqref{eq:soliton} becomes
\[
-\kappa_{\ell\ell}\nu=-6\gamma^\perp
\quad\Longleftrightarrow\quad
\kappa_{\ell\ell}=6\ip{\gamma}{\nu}.
\]
Together with \eqref{eq:lemniscate}, this is the shrinking figure
eight of Edwards et al.\ \cite[Proposition 13]{EdwardsEtAl}.
Ogden and Warren \cite{OgdenWarren2024} also constructed a nontrivial
properly immersed translator for this flow.
\end{remark}

\subsection{Products of Whitney spheres}

\begin{proof}[Proof of Corollary~\ref{cor:product-spheres}]
Write $c_i=2n_i(n_i+2)$ and $f_i=a_i^2f_{W_{n_i}}(u_i,s_i)$.
Equation~\eqref{eq:primitive} gives $F_0^*\lambda_0=d\bigl(\sum_i f_i\bigr)$, so
$F_0$ is exact Lagrangian. Choosing the north or south pole in each factor
gives $2^k$ preimages of the origin, so $F_0$ is not embedded.

The induced metric is $g=\bigoplus_i a_i^2g_i$, where $g_i$ is the metric
\eqref{eq:metric} of $W_{n_i}$. The complex volume form splits, so
$\theta=\sum_i\theta_i$ with $\theta_i$ as in \eqref{eq:angle}. Since
$\Delta_{a^2g_i}=a^{-2}\Delta_{g_i}$, \eqref{eq:lap-angle} gives
\[
\Delta_g\theta+2cf=\sum_{i=1}^k2\Bigl(c-\frac{c_i}{a_i^4}\Bigr)f_i,
\qquad f=\sum_{i=1}^kf_i .
\]
The $f_i$ are nonconstant functions on different factors, so this is
locally constant exactly when \eqref{eq:product-balance} holds, and
Lemma~\ref{lem:criterion} proves the equivalence.
The factor collapse times are $T_i=a_i^4/(4c_i)$, so
$c_i/a_i^4=1/(4T_i)$. Thus \eqref{eq:product-balance} is equivalent to
$T_i=1/(4c)$ for every $i$, and is realized by choosing
$a_i=(c_i/c)^{1/4}$.
\end{proof}

For arbitrary initial scales $a_i$, the velocity in \eqref{eq:flow}
splits over the factors. Applying the scaling formula from
Section~\ref{subsec:whitney} to each factor gives the evolution, up to
reparametrization,
\[
F_t=\bigl(r_1(t)W_{n_1},\ldots,r_k(t)W_{n_k}\bigr),\qquad
r_i(t)=\bigl(a_i^4-8n_i(n_i+2)t\bigr)^{1/4},
\]
for $0\leq t<T_*:=\min_iT_i$.
For any $i$ with $T_i=T_*$, mean curvature scaling gives
\[
\max_M|H_{F_t}|
\geq r_i(t)^{-1}
\max_{\mathbb{S}^{n_i}}|H_{W_{n_i}}|
\longrightarrow\infty
\quad\text{as }t\uparrow T_*.
\]
Since $\max|H_{W_{n_i}}|>0$ by \eqref{eq:angle}, $T_*$ is the first
singular time, with no smooth extension as an immersion of $M$.
If the collapse times differ, the uncollapsed factors define Lagrangian
flows in lower-dimensional ambient spaces between successive collapse
times, with final collapse at $\max_iT_i$.

\subsection{Self-shrinking tori}

\begin{proof}[Proof of Theorem~\ref{thm:circle-tori}]
\proofpart{i}
This follows from Corollary~\ref{cor:product-spheres} with $n_i=a_i=1$
and $c=6$. Writing $F$ for this product immersion, the scaling formula
in Section~\ref{subsec:whitney} gives, for initial immersion $r_0F$ with $r_0>0$,
\[
F_t(\Tn)=\bigl(r_0^4-24t\bigr)^{1/4}F(\Tn),
\qquad 0\leq t<T=\frac{r_0^4}{24},
\]
so the collapse time is independent of $n$ for fixed $r_0$.

\proofpart{ii}
We use the construction $F=z\psi$ in
\cite[(1.2)]{CastroLerma2014}.
Fix an integer $N\geq1$ and coprime positive integers $m_1,m_2$, and set
\begin{equation}\label{eq:nonproduct-curve}
z_N(\sigma)=
\frac{(\sqrt{N+1}+\sqrt N\,e^{i\sigma})^{N+1}}
     {(\sqrt N+\sqrt{N+1}\,e^{i\sigma})^N},
\qquad \sigma\in\R/2\pi\mathbb{Z}.
\end{equation}
The functions $\sqrt{N+1}+\sqrt N\,e^{i\sigma}$ and
$\sqrt N+\sqrt{N+1}\,e^{i\sigma}$ are nonvanishing and have winding
numbers $0$ and $1$, respectively. Hence $z_N$ is a smooth nonzero
closed curve with winding number $-N$ and phase change $-2\pi N$.
Define the closed curve
$\psi_{m_1,m_2}\colon\mathbb{S}^1\to\mathbb{S}^3\subset\C^2$ by
\[
\psi_{m_1,m_2}(\tau)=\left(
\sqrt{\frac{m_2}{m_1+m_2}}\,e^{im_1\tau},\,
\sqrt{\frac{m_1}{m_1+m_2}}\,e^{-im_2\tau}\right),
\qquad \tau\in\R/2\pi\mathbb{Z}.
\]
Define $F_{N;m_1,m_2}\colon\mathbb{T}^2\to\C^2$ by
\begin{equation}\label{eq:weighted-torus}
F_{N;m_1,m_2}(\sigma,\tau)=z_N(\sigma)\psi_{m_1,m_2}(\tau).
\end{equation}
Write $P_N=\sqrt{N(N+1)}$, $A_N=2N+1$ and
$D_N=A_N+2P_N\cos\sigma$.
Since $A_N^2-4P_N^2=1$, we have $D_N>0$. Taking absolute values
and logarithmic derivatives in \eqref{eq:nonproduct-curve} gives
\begin{equation}\label{eq:N-profile-identities}
|z_N|^2=D_N,\qquad
\frac{z_N'}{z_N}=\frac{iP_Ne^{i\sigma}}{D_N},\qquad
|z_N'|^2=\frac{P_N^2}{D_N}.
\end{equation}
In particular, $z_N$ is regular.
For $\psi=\psi_{m_1,m_2}$, we compute
\begin{align*}
&|\psi|=1,\qquad |\psi'|^2=m_1m_2,\qquad
\ip{J\psi}{\psi'}=0,\\
&\det_{\C}(\psi,\psi')
=-i\sqrt{m_1m_2}\,e^{i(m_1-m_2)\tau}.
\end{align*}
Hence $\psi$ is a regular Legendrian curve, and $F_{N;m_1,m_2}$ is a
smooth exact Lagrangian immersion with
\[
g=\frac{P_N^2}{D_N}\,d\sigma^2+m_1m_2D_N\,d\tau^2,\qquad
F_{N;m_1,m_2}^*\lambda_0=df_N,\qquad
f_N=\frac{P_N}{2}\sin\sigma.
\]
Its local Lagrangian angle satisfies
\[
\theta=\arg(z_Nz_N')+(m_1-m_2)\tau-\frac\pi2,\qquad
\theta_\sigma=2-\frac{A_N}{D_N},\quad\theta_\tau=m_1-m_2.
\]
The metric is independent of $\tau$ and $\det g=P_N^2m_1m_2$, so
\[
\Delta_g\theta
=\frac1{P_N^2}\frac{d}{d\sigma}(D_N\theta_\sigma)
=\frac{2D_N'}{P_N^2}
=-\frac4{P_N}\sin\sigma
=-\frac8{P_N^2}f_N.
\]
In particular, $\divg_g(JH)=4\sin\sigma/P_N\not\equiv0$, so these
tori are not Hamiltonian stationary. Lemma~\ref{lem:criterion} gives
\[
J\nabla_g\divg_g(JH)=-c_NF_{N;m_1,m_2}^{\perp},
\qquad c_N=\frac4{P_N^2}=\frac4{N(N+1)}.
\]
The scaling formula in Section~\ref{subsec:whitney} gives the solution,
up to reparametrization,
\[
\left(1-\frac{16t}{N(N+1)}\right)^{1/4}F_{N;m_1,m_2},
\qquad 0\leq t<T_N=\frac{N(N+1)}{16}.
\]

The Gaussian curvature is
\[
K=-\frac{D_N''}{2P_N^2}=\frac{\cos\sigma}{P_N}\not\equiv0,
\]
which excludes a Cartesian product of planar curves.
The torus is not embedded by Theorem~\ref{thm:exactness}. Its image is invariant under
$U_{\delta}=\operatorname{diag}(e^{im_1\delta},e^{-im_2\delta})$, $\delta\in\R$, since
$U_{\delta}F_{N;m_1,m_2}(\sigma,\tau)=F_{N;m_1,m_2}(\sigma,\tau+\delta)$.

To distinguish the parameter $N$, observe that
\[
\frac{\max|F_{N;m_1,m_2}|}{\min|F_{N;m_1,m_2}|}
=\frac{\sqrt{N+1}+\sqrt N}{\sqrt{N+1}-\sqrt N}
=(\sqrt{N+1}+\sqrt N)^2.
\]
This ratio is invariant under unitary transformations, dilations and
reparametrizations, and strictly increases with $N$.
The identity $H=J\nabla_g\theta$ also gives
\[
|H_{N;m_1,m_2}|^2
=\frac{(2D_N-A_N)^2}{P_N^2D_N}
 +\frac{(m_1-m_2)^2}{m_1m_2D_N}.
\]
For $m_1\ne m_2$ this is strictly positive, whereas $H_{N;1,1}=0$
where $D_N=A_N/2$, a value in the range of $D_N$ for every $N\geq1$.
Thus an unequal-weight torus is not congruent to the equal-weight
torus with the same $N$, even after dilation and reparametrization.
For fixed $N$ and $(m_1,m_2)=(1,m)$, the coefficient $(m-1)^2/m$
strictly increases with $m\geq1$, so $\max|H_{N;1,m}|^2$ does as well.
This proves pairwise noncongruence for that subfamily.

For $N=1$, set $z=z_1$ and
$f_{\mathrm{tor}}=f_1=\sin\sigma/\sqrt{2}$; then $c_1=2$ and $T_1=1/8$.
The equal-weight member is unitarily equivalent
to the rotational torus
\[
F(\sigma,\tau)=z(\sigma)(\cos \tau,\sin \tau),
\]
since $F_{1;1,1}=QF$ for the unitary matrix
$Q=\frac1{\sqrt{2}}\left(\begin{smallmatrix}1&i\\1&-i\end{smallmatrix}\right)$.
This immersion is the rotational form of the twisted-torus
construction \cite{ChenCoulibaly}.
In this case the double points can be seen directly from
$z(3\pi/4)=z(5\pi/4)=1$.

\proofpart{iii}
Parametrize $\psi$ by arclength $\tau$, and let $L>0$ denote its length.
The Legendrian condition
$\ip{J\psi}{\psi'}=0$ implies that
the matrix $E(\tau)=(\psi(\tau),\psi'(\tau))$ is unitary. Consequently,
\[
g=|z_\sigma|^2\,d\sigma^2+|z|^2\,d\tau^2,\qquad
F^*\lambda_0=\frac12\operatorname{Im}(\overline z\,z_\sigma)\,d\sigma.
\]
A local Lagrangian angle is
\[
\theta(\sigma,\tau)=\theta_1(\sigma)+\theta_2(\tau),\qquad
\theta_1=\arg(zz_\sigma),\quad \theta_2=\arg\det_{\C}E.
\]
Since $F^*\lambda_0$ has no $d\tau$ component,
Lemma~\ref{lem:criterion} gives $\partial_\tau(\Delta_g\theta)=0$.
On the other hand,
\[
\Delta_g\theta
=\frac1{|z_\sigma||z|}\frac{d}{d\sigma}\left(\frac{|z|}{|z_\sigma|}\,\theta_1'\right)
+\frac{\theta_2''(\tau)}{|z|^2}.
\]
Thus $\theta_2''$ is constant. Periodicity of $\theta_2'$ gives
$\theta_2''=0$ and $\theta_2'=\eta\in\R$.

Differentiating $E^*E=I_2$ and $\det_{\C}E=e^{i\theta_2}$ gives
\[
E'=EB,\qquad
B=\begin{pmatrix}0&-1\\1&i\eta\end{pmatrix},\qquad B^*=-B.
\]
Hence $E(\tau)=E(0)e^{\tau B}$, and $E(\tau+L)=E(\tau)$ gives $e^{LB}=I_2$.
The periodic one-parameter subgroup
\[
U_{\delta}=E(0)e^{\delta B}E(0)^*\in U(2)
\]
therefore satisfies
\[
U_{\delta}\psi(\tau)=\psi(\tau+\delta),\qquad
U_{\delta}F(\sigma,\tau)=F(\sigma,\tau+\delta).
\]
Since $B\ne0$, this subgroup is nontrivial and its image is a circle.
\end{proof}

\begin{remark}[Schoen--Wolfson cones]
Up to the unitary map $\operatorname{diag}(1,i)$, the curve
$\psi_{m_1,m_2}$ is the curve \cite[(5.2)]{SchoenWolfson1999} with
$(p,q)=(m_1,m_2)$, whose cone is Hamiltonian stationary. The torus
$F_{N;m_1,m_2}$ replaces the generating ray of this cone by the closed curve
\eqref{eq:nonproduct-curve}. In part (iii), the equation $E'=EB$ is the
system \cite[(5.1)]{SchoenWolfson1999} with $2a=\eta$, so by the
classification of the Hamiltonian stationary cones
\cite[Section 5]{SchoenWolfson1999} the curve $\psi$ is, after a unitary
transformation and a reparametrization, some $\psi_{p,q}$ with $p,q$
coprime positive integers and $\eta=(p-q)/\sqrt{pq}$.
\end{remark}

\subsection{Exactness of self-shrinkers}

\begin{proof}[Proof of Theorem~\ref{thm:exactness}]
By Lemma~\ref{lem:criterion}, a self-shrinker with constant $c>0$ satisfies
\[
F^*\lambda_0=d\left(\frac{\divg_g(JH)}{2c}\right).
\]
For a shrinker centred at $x_0\in\C^n$, apply the same argument to
$F-x_0$ and use
\[
(F-x_0)^*\lambda_0=F^*\lambda_0-\frac12d\ip{Jx_0}{F}.
\]

If the domain were closed and $F$ embedded, exactness would contradict
Gromov's theorem \cite[Corollary 2.3.B$_2$, p.~330]{Gromov1985}.
\end{proof}

\subsection{Legendrian lifts, integral identities and rigidity}\label{subsec:heisenberg}

For a Legendrian lift $(F,\varphi)$, use the coordinates of
\cite{PigatiRiviere} (defined there for $n=2$):
\[
\zeta=\tfrac12|F|^2+i\varphi,\qquad
\mathfrak r=(2|\zeta|)^{1/2},\qquad
\vartheta=\arg\zeta\in\Bigl[-\frac\pi2,\frac\pi2\Bigr].
\]
Here $\mathfrak r$ is the Kor\'anyi norm of the lift, and
$\vartheta$ is smooth where $\zeta\neq0$; where $F\neq0$ it equals
$\arctan(2\varphi/|F|^2)$, the angular function of \cite{PigatiRiviere}.
\begin{samepage}
\begin{proposition}\label{prop:heisenberg}
Let $F\colon M\to\C^n$ be a self-shrinker with constant $c>0$ on a
closed manifold $M$, with Legendrian lift normalized by
\eqref{eq:legendrian-shrinker}.
\begin{resultparts}
\item
The volume and the height satisfy
\begin{equation}\label{eq:shrinker-volume}
n\operatorname{Vol}(F)=c\int_M\varphi^2\,dV_g .
\end{equation}
\item
If $n=2$, then $\zeta$ has finitely many zeros on $M$, and
\begin{equation}\label{eq:shrinker-angular}
\frac c2\int_M\varphi\,\vartheta\,dV_g
=\int_M|\nabla_g\vartheta|^2\,dV_g+2\pi k,
\end{equation}
where $k$ is the number of zeros of $\zeta$, that is, the number of points
of $M$ which the lift maps to the origin of the Heisenberg group. The
integrand $\varphi\,\vartheta$ is nonnegative.
\end{resultparts}
\end{proposition}
\end{samepage}

\begin{proof}[Proof of Proposition~\ref{prop:heisenberg}]
We first establish the coordinate identities
\begin{equation}\label{eq:zeta}
\Delta_g\zeta=n+i\ip{\nabla_g\theta}{\nabla_g\zeta},
\qquad
|\nabla_g\zeta|^2=2\operatorname{Re}\zeta,
\end{equation}
and, where $\zeta\neq0$,
\begin{equation}\label{eq:log-zeta}
\Delta_g\log|\zeta|
=\frac{n-2}{2}\bigl|\nabla_g\log|\zeta|\bigr|^2
+\frac{n+2}{2}|\nabla_g\vartheta|^2
-\ip{\nabla_g\theta}{\nabla_g\vartheta}.
\end{equation}
These identities hold for any Legendrian lift, without the shrinker
assumption, and will also be used in the proof of
Theorem~\ref{thm:koranyi-rigidity}. Here $\ip{\cdot}{\cdot}$ is extended
complex bilinearly, whereas
$|\nabla_g\zeta|^2
=|\nabla_g\operatorname{Re}\zeta|^2+|\nabla_g\operatorname{Im}\zeta|^2$.
We have $\nabla_g\bigl(\frac12|F|^2\bigr)=F^\top$, and
$d\varphi=2F^*\lambda_0=\ip{JF}{dF}$ gives
$\nabla_g\varphi=(JF)^\top=J(F^\perp)$. Hence
$|\nabla_g\zeta|^2=|F^\top|^2+|F^\perp|^2=|F|^2=2\operatorname{Re}\zeta$.
Taking divergences, using $\Delta_gF=H=J\nabla_g\theta$ and
$\ip{Je_j}{e_j}=0$ for an orthonormal tangent frame, gives
\[
\Delta_g\bigl(\tfrac12|F|^2\bigr)=n-\ip{\nabla_g\theta}{\nabla_g\varphi},
\qquad
\Delta_g\varphi=\ip{\nabla_g\theta}{\nabla_g(\tfrac12|F|^2)},
\]
which is the first identity in \eqref{eq:zeta}. Where $\zeta\neq0$, a local
branch $\log\zeta=\log|\zeta|+i\vartheta$ satisfies
\[
\Delta_g\log\zeta
=\frac{\Delta_g\zeta}{\zeta}-\frac{\ip{\nabla_g\zeta}{\nabla_g\zeta}}{\zeta^2}
=\frac n\zeta+i\ip{\nabla_g\theta}{\nabla_g\log\zeta}
-\ip{\nabla_g\log\zeta}{\nabla_g\log\zeta},
\]
and, by the second identity in \eqref{eq:zeta},
\[
\operatorname{Re}\frac1\zeta=\frac{\operatorname{Re}\zeta}{|\zeta|^2}
=\frac{|\nabla_g\zeta|^2}{2|\zeta|^2}
=\frac12\Bigl(\bigl|\nabla_g\log|\zeta|\bigr|^2+|\nabla_g\vartheta|^2\Bigr).
\]
Taking real parts gives \eqref{eq:log-zeta}.

\proofpart{i}
Integrating the real part of the first identity in \eqref{eq:zeta} over
$M$ and using \eqref{eq:legendrian-shrinker},
\[
n\operatorname{Vol}(F)=\int_M\ip{\nabla_g\theta}{\nabla_g\varphi}\,dV_g
=-\int_M\varphi\,\Delta_g\theta\,dV_g=c\int_M\varphi^2\,dV_g .
\]

\proofpart{ii}
Let $\zeta(x_0)=0$, so that $F(x_0)=0$ and $\varphi(x_0)=0$, and let $x$ be
normal coordinates centred at $x_0$. Then $|F|^2=|x|^2+O(|x|^3)$. Since
$F=dF_{x_0}(x)+O(|x|^2)$ and the image of $dF_{x_0}$ is a Lagrangian
plane, $d\varphi=\ip{JF}{dF}=O(|x|^2)$, so $\varphi=O(|x|^3)$. Hence
$x_0$ is an isolated zero of $\zeta$, the function $\vartheta$ is
Lipschitz near $x_0$, and $\log|\zeta|-2\log|x|$ has bounded gradient near
$x_0$. In particular, $\zeta$ has finitely many zeros.

For $n=2$, \eqref{eq:log-zeta} reads
$\Delta_g\log|\zeta|=2|\nabla_g\vartheta|^2
-\ip{\nabla_g\theta}{\nabla_g\vartheta}$.
We integrate this over the complement of the geodesic discs of radius
$\varepsilon$ about the zeros of $\zeta$ and let $\varepsilon\to0$. By the
divergence theorem and the expansion above, the left-hand side tends to
$-4\pi k$. Since $\vartheta$ is bounded and Lipschitz near the zeros,
\[
\int_M\ip{\nabla_g\theta}{\nabla_g\vartheta}\,dV_g
=-\int_M\vartheta\,\Delta_g\theta\,dV_g
=c\int_M\varphi\,\vartheta\,dV_g
\]
by \eqref{eq:legendrian-shrinker}. This gives \eqref{eq:shrinker-angular}. Finally,
$\vartheta=\arg\zeta$ has the sign of $\varphi=\operatorname{Im}\zeta$.
\end{proof}

\begin{proof}[Proof of Theorem~\ref{thm:koranyi-rigidity}]
The condition $(F,\varphi)(M)\subset S_K(R)$ gives
\[
\zeta=\frac{R^2}{2}e^{i\vartheta},\qquad
|F|^2=R^2\cos\vartheta,\qquad
\varphi=\frac{R^2}{2}\sin\vartheta.
\]
Since $\zeta$ never vanishes and has nonnegative real part,
$\vartheta\colon M\to[-\pi/2,\pi/2]$ is globally smooth.
The second identity in \eqref{eq:zeta} gives
\begin{equation}\label{eq:constant-gauge-angle}
|\nabla_g\vartheta|^2=\frac4{R^2}\cos\vartheta.
\end{equation}
The function $\vartheta$ is nonconstant, since otherwise this identity
would force $F\equiv0$. Its extrema are therefore $-\pi/2$ and
$\pi/2$. Its critical set is $Z=F^{-1}(0)$, which is finite: the
immersion $F$ is locally an embedding, so $Z$ is discrete, and it is
closed in the compact manifold $M$.

On $M_*=M\setminus Z$, put
$X=\nabla_g\vartheta/|\nabla_g\vartheta|^2$, so $X\vartheta=1$.
With $V=dF(\nabla_g\vartheta)$, the identities
$\nabla_g(\frac12|F|^2)=F^\top$ and
$\nabla_g\varphi=J(F^\perp)$ give
\[
F^\top=-\frac{R^2}{2}\sin\vartheta\,V,\qquad
F^\perp=-\frac{R^2}{2}\cos\vartheta\,JV.
\]
Consequently,
\[
dF(X)=\left(-\frac12\tan\vartheta+\frac i2\right)F.
\]
The flow of $X$ exists between any two levels in
$(-\pi/2,\pi/2)$, since their intervening level strip is compact.
It identifies $M_*$ with $(-\pi/2,\pi/2)\times\Sigma$, where
$\Sigma=\{\vartheta=0\}$. Integrating the preceding equation gives
\begin{equation}\label{eq:constant-gauge-profile}
F(\vartheta,q)=R\sqrt{\cos\vartheta}\,e^{i\vartheta/2}\psi(q),
\qquad \psi(q)=F(0,q)/R\in\mathbb{S}^{2n-1}.
\end{equation}

By \eqref{eq:constant-gauge-angle},
$|X|=R/(2\sqrt{\cos\vartheta})$, which is integrable at both endpoints.
Thus the flow lines converge uniformly in $M$ as $\vartheta\to\pm\pi/2$,
defining continuous endpoint maps $\Sigma\to Z$.
Fix a connected component $\Sigma_0$ of $\Sigma$. Its lower endpoint
map is constant, since $Z$ is finite; denote this endpoint by $p$.
Smoothness of the immersion at $p$ implies that every limiting unit
direction of $F$ lies in the Lagrangian plane $dF_p(T_pM)$.
Equation \eqref{eq:constant-gauge-profile} therefore gives
\[
e^{-i\pi/4}\psi(q)\in dF_p(T_pM),\qquad q\in\Sigma_0.
\]
In particular, $\psi(\Sigma_0)$ lies in the fixed Lagrangian plane
$P_{\Sigma_0}=e^{i\pi/4}dF_p(T_pM)$.

Define $s=\tan(\vartheta/2)$ and $Y=(1-is)F/R$ on all of $M$.
These smooth functions satisfy
\[
|Y|^2+s^2=1,\qquad
Y(\vartheta,q)=\sqrt{1-s^2}\,\psi(q)\quad\text{on }M_*.
\]
Hence $Y$ maps each component of $M_*$ into a fixed Lagrangian
plane. These planes agree across $Z$. Indeed, at any $p\in Z$,
\[
dY_p=\frac{1-is(p)}R\,dF_p
\]
has rank $n$. Taking limits of $dY$ from an adjacent component shows
that its plane must equal $dY_p(T_pM)$. Since $Z$ is finite, the
planes thus extend to a locally constant choice on $M$.
Connectedness gives a single Lagrangian plane $P$ containing $Y(M)$.

The map $G=(Y,s)$ takes values in the unit sphere
$\mathbb{S}(P\times\R)$. It is an immersion, because
\[
F=R\,\frac{1+is}{1+s^2}\,Y
\]
and $dG(v)=0$ would imply $dF(v)=0$.
Choose $U\in U(n)$ with $U(\R^n)=P$ and set
$\Pi=(U^{-1}Y,s)\colon M\to\mathbb{S}^n$. Then $\Pi$ is a local
diffeomorphism, and the displayed formula gives
$F=R\,U\circ W_n\circ\Pi$.

By compactness, $\Pi\colon M\to\mathbb{S}^n$ is a finite smooth covering.
For $n\geq2$, connectedness of $M$ and simple connectedness of
$\mathbb{S}^n$ imply that $\Pi$ is a diffeomorphism. For $n=1$, $F$ is
a finite cover of a scaled and rotated Bernoulli lemniscate.

In both cases, the shrinker equation is preserved by unitary
transformations and local reparametrizations.
Theorem~\ref{thm:main} and scaling therefore give $c=2n(n+2)/R^4$.
\end{proof}

For the Whitney immersion $W_n$, equations \eqref{eq:whitney},
\eqref{eq:angle} and \eqref{eq:primitive} give
\[
\varphi=\frac{s}{1+s^2},\qquad
\zeta=\frac12\,\frac{1+is}{1-is}=\tfrac12e^{i\vartheta},\qquad
\vartheta=2\arctan s,\qquad
\theta=\frac{n+2}{2}\,\vartheta+\frac\pi2 .
\]
Thus the $\zeta$-image is the right semicircle of radius $1/2$, and the
lift lies in the Kor\'anyi sphere $\{\mathfrak r=1\}$. The angle relation
makes the last two terms of \eqref{eq:log-zeta} cancel pointwise.
For $W_2$, we have $k=0$ and both sides of
\eqref{eq:shrinker-angular} equal $8\pi$.

For the tori $F_{N;m_1,m_2}$, the formulas in the proof of
Theorem~\ref{thm:circle-tori}(ii) give
\[
\varphi_N=P_N\sin\sigma,\qquad
\zeta_N=\frac12D_N+iP_N\sin\sigma
=\frac{2N+1}{2}+\sqrt{N(N+1)}\,e^{i\sigma}.
\]
Since $(2N+1)^2/4-N(N+1)=1/4$, each $\zeta_N$-image is a circle
lying strictly in the right half-plane. In particular, $k=0$ for
every member of this family. Moreover,
$dV_g=P_N\sqrt{m_1m_2}\,d\sigma\,d\tau$ gives
\[
\int_{\mathbb{T}^2}\varphi_N^2\,dV_g
=\frac{P_N^2}{2}\operatorname{Vol}(F_{N;m_1,m_2}),
\]
which verifies \eqref{eq:shrinker-volume} with $n=2$ and
$c_N=4/P_N^2$. For the product of two lemniscates in
Theorem~\ref{thm:circle-tori}(i), exactly two of the four preimages of
$F=0$ also have $\varphi=0$, so $k=2$ in \eqref{eq:shrinker-angular}.

\begin{remark}[Relation to Hamiltonian stationary monotonicity]
In dimension two, the coordinates above and the identity
\eqref{eq:log-zeta} are related to the monotonicity formulas for
Hamiltonian stationary Legendrian surfaces
\cite{SchoenWolfson1999,SchoenWolfson2001}; see also
\cite[Proposition 3.4 and Theorem 3.6]{PigatiRiviere}.
The mixed term in \eqref{eq:log-zeta} is
$-\ip{\nabla_g\theta}{\nabla_g\vartheta}$. On $\{\zeta\neq0\}$,
\[
-\ip{\nabla_g\theta}{\nabla_g\vartheta}
=-\divg_g(\vartheta\nabla_g\theta)+\vartheta\Delta_g\theta.
\]
For a Hamiltonian stationary immersion, $\Delta_g\theta=0$, so
this term equals $-\divg_g(\vartheta\nabla_g\theta)$. For a
self-shrinker, integration by parts using $\Delta_g\theta=-c\varphi$
produces the left-hand side of \eqref{eq:shrinker-angular}.
\end{remark}

\section{Exact deformations and Hamiltonian isotopy}\label{sec:hamiltonian}

An \emph{ambient Hamiltonian isotopy} is a family of symplectic
diffeomorphisms $\Phi_t\colon\C^2\to\C^2$, with
$\Phi_0=\operatorname{id}$, generated by a time-dependent
Hamiltonian vector field. Such an isotopy preserves the pairs of
coincident points: for $G_t=\Phi_t\circ G_0$,
\[
G_t(x)=G_t(y)\quad\Longleftrightarrow\quad G_0(x)=G_0(y).
\]

A smooth family $G_t\colon M\to\C^2$ of Lagrangian immersions is
an \emph{exact Lagrangian deformation} if the one-form
\[
\beta_t(X)=\omega_0\bigl(\partial_tG_t,dG_t(X)\bigr)
\]
is exact on $M$ for every $t$. The flow \eqref{eq:flow} has
$\beta_t=-d\bigl(\divg_{g(t)}(JH_t)\bigr)$.
Here the Hamiltonian potentials are functions on $M$. An exact
Lagrangian deformation need not arise from an ambient Hamiltonian
isotopy; see \cite[Section 2.2]{ChenWarren}.

\begin{proposition}\label{prop:deformation-obstructions}
Let $F=F_{N;m_1,m_2}$ be a torus in \eqref{eq:weighted-torus}.
All comparisons below allow reparametrization of the domain.
\begin{resultparts}
\item
With the convention $\mu_F=[d\theta/\pi]$, its Maslov indices on
the positively oriented coordinate circles are
\[
\bigl(\mu_F([S^1_\sigma]),\mu_F([S^1_\tau])\bigr)
=\bigl(2(1-2N),\,2(m_1-m_2)\bigr).
\]
The positive generator of $\mu_F(H_1(\mathbb{T}^2;\mathbb{Z}))$ is
\[
d_F=2\gcd(2N-1,|m_1-m_2|).
\]
Two members with different $d_F$ cannot be joined by a smooth family
of Lagrangian immersions. In particular, the equal-weight tori
$F_{N;1,1}$ for different $N$ lie in distinct deformation classes.
\item
The torus $F$ cannot be joined to an embedded Lagrangian torus by an
exact Lagrangian deformation, or to any torus with vanishing Maslov
class by a smooth family of Lagrangian immersions. The latter
excludes the product of two Bernoulli lemniscates.
\item
Positive dilations of $F$ give exact Lagrangian deformations, but
$F$ and $rF$ are not ambient Hamiltonian isotopic for $r>0$, $r\ne1$.
\end{resultparts}
\end{proposition}

\begin{proof}[Proof of Proposition~\ref{prop:deformation-obstructions}]
\proofpart{i}
Recall $A_N=2N+1$, $P_N=\sqrt{N(N+1)}$ and
$D_N=A_N+2P_N\cos\sigma$ as in the proof of
Theorem~\ref{thm:circle-tori}(ii). Since $A_N^2-4P_N^2=1$,
\[
\int_0^{2\pi}\frac{d\sigma}{D_N}
=\frac{2\pi}{\sqrt{A_N^2-4P_N^2}}=2\pi.
\]
The angle derivatives $\theta_\sigma=2-A_N/D_N$ and
$\theta_\tau=m_1-m_2$ therefore give
\[
\frac1\pi\int_0^{2\pi}\theta_\sigma\,d\sigma
=4-2A_N=2(1-2N),\qquad
\frac1\pi\int_0^{2\pi}\theta_\tau\,d\tau=2(m_1-m_2).
\]
These identities give the stated Maslov indices and $d_F$.
In a smooth family of Lagrangian immersions, the periods of
$d\theta/\pi$ vary continuously and are integer-valued, so the
Maslov class is constant. Reparametrization acts by an automorphism
of $H_1(\mathbb{T}^2;\mathbb{Z})$ and preserves the image of $\mu_F$,
hence $d_F$. For $m_1=m_2=1$, this invariant is $2(2N-1)$ and
strictly increases with $N$.

\proofpart{ii}
Let $G_t$ be an exact Lagrangian deformation with $G_0=F$.
Cartan's formula gives
\[
\frac{d}{dt}(G_t^*\lambda_0)
=\beta_t+d\bigl(\lambda_0(\partial_tG_t)\bigr),
\]
so $[G_t^*\lambda_0]\in H^1(\mathbb{T}^2;\R)$ is independent of $t$.
Since $F$ is exact, every $G_t$ is exact. Gromov's theorem
\cite[Corollary 2.3.B$_2$]{Gromov1985}
therefore rules out every embedded endpoint, including products of
round circles and Chekanov tori.

Part (i) gives $\mu_F\ne0$ for every $N\geq1$, whereas the product of
two Bernoulli lemniscates has a globally defined real-valued angle by
\eqref{eq:angle}. Invariance of the Maslov class proves the remaining
assertions.

\proofpart{iii}
Recall that $F^*\lambda_0=df_N$, with $f_N=(P_N/2)\sin\sigma$.
For $G_t=r(t)F$ with $r(t)>0$, we have
$\beta_t=2r(t)r'(t)\,df_N$, so dilations give exact Lagrangian
deformations. To distinguish $F$ from its dilations under ambient
Hamiltonian isotopy, we exhibit a self-intersection with nonzero
primitive difference.
Let $\alpha_N$ be the continuous argument of $z_N$ on $[0,\pi]$
with $\alpha_N(0)=0$. By \eqref{eq:N-profile-identities},
\[
\alpha_N'=\frac{P_N\cos\sigma}{D_N}
=\frac12\left(1-\frac{A_N}{D_N}\right),\qquad
\alpha_N(\pi)=-N\pi.
\]
Here the endpoint value follows by integrating, since
$\int_0^\pi D_N^{-1}\,d\sigma=\pi$.
The function $\alpha_N$ is strictly increasing on $(0,\pi/2)$ and
strictly decreasing on $(\pi/2,\pi)$, so it has a unique zero
$\sigma_N\in(\pi/2,\pi)$. Since
$z_N(2\pi-\sigma)=\overline{z_N(\sigma)}$, for every $\tau$ the points
$x=(\sigma_N,\tau)$ and $y=(2\pi-\sigma_N,\tau)$ satisfy
\[
F(x)=F(y),\qquad f_N(x)-f_N(y)=P_N\sin\sigma_N>0.
\]
Consequently the quantity
\[
\mathcal A(F)=\max\{\,|f_N(x)-f_N(y)|:F(x)=F(y)\,\}
\]
exists by compactness and is positive.
An ambient Hamiltonian diffeomorphism $\Phi$ satisfies
$\Phi^*\lambda_0=\lambda_0+dS$, so the primitive of $\Phi\circ F$ is
$f_N+S\circ F$. The map $\Phi$ preserves the coincident pairs, and
$S\circ F$ has equal values at each such pair. Hence $\mathcal A$ is
invariant under ambient Hamiltonian diffeomorphisms. It is also
unchanged by reparametrization, while $\mathcal A(rF)=r^2\mathcal A(F)$.
Since $\mathcal A(F)>0$, an ambient Hamiltonian isotopy from $F$ to $rF$
is impossible for $r\ne1$.
\end{proof}

\end{document}